\documentclass{amsart}
\usepackage{amsmath}
\usepackage{amsthm}
\usepackage{amsfonts}
\usepackage{amscd}
\usepackage{mathrsfs}
\usepackage{amssymb}
\usepackage{euscript}
\usepackage{oldgerm}
\usepackage{graphicx}
\frakfamily

\theoremstyle{plain}
\newtheorem{thm}{Theorem}[section]
\newtheorem{prop}[thm]{Proposition}
\newtheorem{cor}[thm]{Corollary}
\newtheorem{lem}[thm]{Lemma}

\theoremstyle{definition}
\newtheorem{notation}[thm]{Notation}
\newtheorem{defn}[thm]{Definition}
\theoremstyle{remark}
\newtheorem{remark}[thm]{Remark}

\newtheorem{example}[thm]{Example}
\renewenvironment{proof}[1][Proof]{\textbf{#1.} }{\hfill \rule{0.5em}{0.5em}}
\numberwithin{equation}{section}

\title{Malliavin calculus for Lie group-valued Wiener functions}
\author{Tai Melcher}
\address{Department of Mathematics, University
of Virginia, Charlottesville, VA 22936}
\thanks{The author was supported in part by NSF Grants 99-71036 and DMS
0202939.} 
\email{melcher@virginia.edu}


\begin{document}

\begin{abstract}
Let $G$ be a Lie group equipped with a set of left invariant vector fields.  These vector fields generate a function $\xi$ on Wiener space into $G$  via the stochastic version of Cartan's rolling map.  It is shown here that, for any smooth function $f$ with compact support, $f(\xi)$ is Malliavin differentiable to all orders and these derivatives belong to $L^p(\mu)$ for all $p>1$, where $\mu$ is Wiener measure.
\end{abstract}

\maketitle
\setcounter{tocdepth}{1}
\tableofcontents

\section{Introduction}
\label{s:intro}

Malliavin \cite{Malliavin78b,Malliavin78c} first introduced the notion of derivatives of Wiener functionals as part of a program for producing a probabilistic proof of the celebrated H\"ormander Theorem, which states that solutions to certain stochastic differential equations have smooth transition densities.  The purpose of the present paper is to show that, under certain conditions, functions of these solutions are ``Malliavin smooth," in the sense that they belong to Sobolev spaces of all orders.   Malliavin calculus has found applications in many aspects of classical and stochastic analysis.  In particular, applications of the following result to heat kernel inequalities of hypoelliptic operators can be found in \cite{Melcher2}.

Let $G$ be a Lie group with identity $e$, and let $\{\tilde{X}_i\}_{i=1}^k$ be a set of left invariant vector fields.  Let $\mathscr{W}(\mathbb{R}^k)$ denote standard $k$-dimensional path space equipped with Wiener measure $\mu$,  and let $\xi:[0,1]\times \mathscr{W}(\mathbb{R}^k) \rightarrow G$ denote the solution to the Stratonovich stochastic differential equation
\[ d\xi_t = \sum_{i=1}^k \tilde{X}_i(\xi_t)\circ db_t^i \text{ with } \xi_0=e, \]
where $(b^1,\ldots,b^k)$ is a $k$-dimensional Brownian motion.  Then, for all $f\in C_c^\infty(G)$, $f(\xi_t)$ is Malliavin differentiable to all orders, and these derivatives are in $L^p(\mu)$ for all $p>1$.  Theorem \ref{t:xi.smooth} states this result explicitly and gives an expression for the first order derivative of $f(\xi_t)$.

Regularity results of this type have been known for manifolds equipped with vector fields satisfying certain boundedness conditions.  For example, $f(\xi_t)$ is Malliavin smooth when $M$ is a compact Riemannian manifold, or, more generally, when $M$ is diffeomorphic to a closed submanifold of $\mathbb{R}^N$, and this diffeomorphism induces a mapping on the vector fields to vector fields on $\mathbb{R}^N$ which have at most linear growth and bounded derivatives of all orders; see Taniguchi \cite{Taniguchi83}.  These previous results are not sufficient for a set of general left invariant vector fields on a Lie group, as the following example demonstrates.

\begin{example}
Let $\mathfrak{g}$ be the $4$-dimensional algebra defined in a basis $\{X_1,X_2,X_3,X_4\}$ by the brackets
\[ [X_1,X_2]=X_3 \qquad \text{and} \qquad [X_1,X_3]=X_4, \]
where the other brackets (except those obtained by anticommutativity) are 0.  It can be shown (for example, via the Baker-Campbell-Hausdorff formula) that the associated Lie group may be realized as $\mathbb{R}^4$ with the group operation given componentwise as
\begin{equation*}
\begin{split}
[(a,b,c,d)\cdot(a',b',c',d')]_1 &= a+a' \\
[(a,b,c,d)\cdot(a',b',c',d')]_2 &= b+b' \\ 
[(a,b,c,d)\cdot(a',b',c',d')]_3 &= c+c' + \frac{1}{2}(ab'-a'b) \\ 
[(a,b,c,d)\cdot(a',b',c',d')]_4 &= d+d' + \frac{1}{2}(ac'-a'c) \\ 
	&\qquad  + \frac{1}{12}(a^2b'-aa'b - aa'b' + (a')^2b)
\end{split}
\end{equation*}
Now let $X=(1,0,0,0)$  at $e=(0,0,0,0)$ in $G$.  Extending $X$ to a left invariant vector field in the usual way gives,
\begin{equation*}
\begin{split}
X(a,b,c,d) &= L_{(a,b,c,d)*}X_1
	= \frac{d}{d\varepsilon}\bigg|_{\varepsilon=0} (a,b,c,d)\cdot(\varepsilon,0,0,0) \\
	&=  \frac{d}{d\varepsilon}\bigg|_{\varepsilon=0} (a+\varepsilon, b, c - \frac{1}{2}b\varepsilon, 
		d - \frac{1}{2}c\varepsilon + \frac{1}{12}(-ab\varepsilon+b\varepsilon^2)) \\
	&= (1,0,-\frac{1}{2}b, -\frac{1}{2}c - \frac{1}{12}ab).
\end{split}
\end{equation*}
So if $(w,x,y,z)$ are the standard coordinates in $G=\mathbb{R}^4$,
\[ X = \frac{\partial}{\partial w} - \frac{1}{2}x\frac{\partial}{\partial y} - \left( \frac{1}{2}y + \frac{1}{12}wx\right) \frac{\partial}{\partial z}. \]
Thus, the standard coordinate representation for this vector field has non-linear coefficients.
\end{example}
 
Section \ref{s:Mall.def} sets up the notation for this paper and recalls standard definitions and notions of differentiability on Wiener space.  Section \ref{s:main} gives the main result in Theorem \ref{t:xi.smooth}, and Section \ref{s:SDEp} presents the proof.

\subsection*{Acknowledgement.} I thank Bruce Driver for suggesting
this problem and for many valuable discussions throughout the
preparation of this work.

\section{Background and main result}

\subsection{Wiener space calculus}
\label{s:Mall.def}
This section contains a brief introduction to basic Wiener space definitions of differentiability.  For a more complete exposition, consult \cite{Ikeda84,IkedaWatanabe89,KS1,KS2,Malliavin78,Malliavin97,Norris86b,Nualart95,Watanabe84,Watanabe97} and references contained therein.  In particular, the first two chapters of Nualart \cite{Nualart95} and Chapter V of Ikeda and Watanabe \cite{IkedaWatanabe89} are cited often here.

Let $(\mathscr{W}(\mathbb{R}^{k}),\mathcal{F},\mu)$ denote classical k-dimensional
Wiener space. That is, $\mathscr{W}=\mathscr{W}(\mathbb{R}^k)$ is the Banach space of continuous paths $\omega:[0,1]\rightarrow\mathbb{R}^{k}$ such that $\omega_0=0$, equipped with
the supremum norm
\[
\Vert\omega\Vert=\max_{t\in\lbrack0,1]}|\omega_t|,
\]
$\mu$ is standard Wiener measure, and $\mathcal{F}$ is the completion of the
Borel $\sigma$-field on $\mathscr{W}$ with respect to $\mu$.  By definition of $\mu,$ the process
\[
b_t(\omega) = (b_t^1(\omega), \ldots , b_t^k(\omega) ) = \omega_t
\]
is an $\mathbb{R}^k$ Brownian motion. For those $\omega\in \mathscr{W}$ which are
absolutely continuous, let
\[
E(\omega):=\int_{0}^{1}|\dot{\omega}_s|^{2}\,ds
\]
denote the energy of $\omega$. The Cameron-Martin space is the Hilbert space
of finite energy paths,
\[
\mathscr{H}=\mathscr{H}(\mathbb{R}^k):=\{\omega\in \mathscr{W}(\mathbb{R}^k):\omega\text{ is
absolutely continuous and }E(\omega)<\infty\},
\]
equipped with the inner product
\[
(h,k)_{\mathscr{H}}:=\int_{0}^{1} \dot{h}_s\cdot\dot{k}_s\,ds, \quad\text{ for all }h,k\in
\mathscr{H}.
\]
More generally, for any finite dimensional vector space $V$ equipped with an inner product, let $\mathscr{W}(V)$  denote path space on $V$, and $\mathscr{H}(V)$ denote the set of Cameron-Martin paths, where the definitions are completely analogous, replacing the inner products and norms where necessary.
%We may identify the Cameron-Martin space $H^1$ with $H=L^{2}([0,1],\mathbb{R}^k)$
%in the obvious way
%\[
%h\in H^{1} \mapsto\dot{h}\in H.
%\]
%With the standard inner product on $H$,
%\[ (h,k)_H = \int_0^1 h(s)\cdot k(s)\,ds, \quad\forall~h,k\in H, \]
%this mapping is an isometry, and the spaces are isomorphic.  In the sequel, we make this
%identification without further comment.
\iffalse
To define a notion of differentiation for functions on $\mathscr{W}$, let $B=\{B(h),h\in H\}$ be the Gaussian random field given by the It\^o integral
\[
B(h)=\int_{0}^{1}\dot{h}(t)\cdot db_{t}.
\]
$B$ is often called the Wiener integral and is well-defined with $B(h)\in L^2(\mu)$.
%$B$ is an isonormal Gaussian process associated to the Hilbert space H.
\fi

\begin{defn}
Denote by $\mathcal{S}$ the class of {\em smooth cylinder functionals}, 
random variables $F:\mathscr{W}\rightarrow\mathbb{R}$ such that
\begin{equation}
\label{e:cyl}
F(\omega)=f(\omega_{t_1},\ldots,\omega_{t_n}),
\end{equation}
for some $n\geq1,$ $0 < t_{1} < \cdots < t_{n} \le 1$, and function $f\in
C^\infty_{p}((\mathbb{R}^{k})^n)$, the functions $f\in C^\infty(\mathbb{R}^{kn})$ such that $f$ and all of its partial derivatives have at most polynomial growth.  For $E$ be a real separable Hilbert space, let $\mathcal{S}_{E}$ be the set of $E$-valued smooth cylinder functions $F:\mathscr{W}\rightarrow E$ of the form
\begin{equation}
\label{e:Ecyl}
F=\sum_{j=1}^{m} F_{j} e_{j},
\end{equation}
for some $m\ge1$, $e_{j}\in E$, and $F_{j}\in\mathcal{S}$.
\end{defn}

\begin{defn}
Fix $h\in\mathscr{H}$.  The directional derivative of a smooth cylinder functional $F\in\mathcal{S}$ of the form (\ref{e:cyl}) along $h$ is given by
\begin{equation*}
\partial_{h}F(\omega):=\frac{d}{d\epsilon}\bigg|_{0}F(\omega+\epsilon h) = 
  \sum_{i=1}^n  \nabla^i f(\omega_{t_1},\ldots,\omega_{t_n})\cdot h_{t_i},
\end{equation*}
where $\nabla^i f$ is the gradient of $f$ with respect to the $i^{\mathrm{th}}$ variable.
\end{defn}

The following integration by parts result is standard; see for example Theorem 8.2.2 of Hsu \cite{Hsu03}.
\begin{prop}
\label{p:dh*}
Let $F,G\in\mathcal{S}$ and $h\in\mathscr{H}$.  Then
\[ (\partial_h F,G)_\mathscr{H} = (F,\partial_h^* G)_\mathscr{H}, \]
where $\partial_h^*=-\partial_h + \int_0^1 \dot{h}_s\cdot db_s$.
\end{prop}

\begin{defn}
The gradient of a smooth cylinder functional $F\in\mathcal{S}$ is the random process $D_tF$
taking values in $\mathscr{H}$ such that $(DF,h)_{\mathscr{H}}=\partial_{h}F$.  It may be determined that, for $F$ of the form (\ref{e:cyl}),
\[
D_{t}F = \sum_{i=1}^{n} \nabla^if(\omega_{t_1},\ldots,\omega_{t_n}) (t_i\wedge t),
\]
where $s\wedge t=\min\{s,t\}$.  For $F\in\mathcal{S}_E$ of the form (\ref{e:Ecyl}), define the derivative $D_tF$ to be the random process taking values in $\mathscr{H}\otimes E$ given by
\[ D_tF := \sum_{j=1}^m D_tF\otimes e_j. \]
\end{defn}

Iterations of the derivative for smooth functionals $F\in\mathcal{S}$ are given by
\[
D_{t_{1},\ldots,t_{k}}^{k}F=D_{t_{1}}\cdots D_{t_{k}}F \in \mathscr{H}^{\otimes k},
\]
for $k\in\mathbb{N}$.  For $F\in\mathcal{S}_E$,
\[
D^{k} F = \sum_{j=1}^{m} D^{k} F_{j} \otimes e_{j},
\]
and these are measurable functions defined almost everywhere on $[0,1]^{k}\times \mathscr{W}$.  The operator $D$ on $\mathcal{S}_E$ is closable, and there exist closed extensions $D^{k}$ to $L^{p}(\mathscr{W},\mathscr{H}^{\otimes k}\otimes E)$; see, for example \cite{Nualart95}, Theorem 8.28 of \cite{Hsu03}, or Theorem 8.5 of \cite{IkedaWatanabe89}.
% Nualart: see page 27, item (ii) 
Denote the closure of the derivative operator also by $D$ and the domain of $D^{k}$ in $L^{p}([0,1]^{k}\times \mathscr{W})$ by ${\mathcal{D}}^{k,p}$, which is the completion of the family of smooth Wiener functionals $\mathcal{S}$ with respect to the seminorm $\Vert\cdot\Vert_{k,p,E}$ on $\mathcal{S}_E$ given by
\[
\|F\|_{k,p,E} := \left( \sum_{j=0}^{k} \mathbb{E}(\|D^{j}F\|
   _{\mathscr{H}^{\otimes j}\otimes E}^{p}) \right)^{1/p},
\]
for any $p\ge1$.  Let
\[
\mathcal{D}^{k,\infty}(E) := \bigcap_{p>1} \mathcal{D}^{k,p}(E) \text{ and } 
   \mathcal{D}^{\infty}(E) := \bigcap_{p>1}\bigcap_{k\ge1}{\mathcal{D}}^{k,p}(E).
\]
When $E=\mathbb{R}$, write $\mathcal{D}^{k,p}(\mathbb{R})=\mathcal{D}^{k,p}$, $\mathcal{D}^{k,\infty}(\mathbb{R})=\mathcal{D}^{k,\infty}$, and $\mathcal{D}^\infty(\mathbb{R})=\mathcal{D}^\infty$.  Also let  
\[ L^{\infty-}(\mu) := \bigcap_{p>1}L^{p}(\mu). \]
The notion of Sobolev spaces of Wiener functionals was first introduced by Shigekawa \cite{Shigekawa80} and Stroock \cite{Stroock80a,Stroock80b}.

The operator $\partial_h$ on $\mathcal{S}$ is also closable, and there exists a closed extension of $\partial_{h^1}\cdots\partial_{h^k}$ to $L^p(\mu)$.  Denote the closure of $\partial_h$ also by $\partial_h$, with domain $\mathrm{Dom}(\partial_h)$.  Denote by $G^{k,p}$ the class of functions $F\in L^p(\mu)$ such that $\partial_{h^1}\cdots\partial_{h^j}F\in L^p(\mu)$, for all $h^1,\ldots,h^j\in \mathscr{H}$, $j=1,\ldots, k$.  The norm on $G^{k,p}$ is  given by
\[ \|F\|_{G^{k,p}} := \sum_{j=0}^k \left(\mathbb{E}\|D^jF \|_{(\mathscr{H}^{\otimes j})^*}^p\right)^{1/p}, \]
where 
\[ \|D^jF\|_{(\mathscr{H}^{\otimes j})^*} = \sup\{|\partial_{h^1}\cdots\partial_{h^j}F |: h^i\in \mathscr{H}, |h^i|_\mathscr{H}\le1 \}  \]
is the operator norm on the space of continuous linear functionals $(\mathscr{H}^{\otimes j})^*$.

The following result follows from Proposition 5.4.6 and Corollary 5.4.7 of \cite{Bogachev98}.
(Note that the space $\mathcal{D}^{k,p}$ here corresponds to the space $W^{k,p}$ in that text.)

\begin{thm}
\label{t:diff}
For all $k\in\mathbb{N}$ and $p\in(1,\infty)$, $\mathcal{D}^{k,p}\subset G^{k,p}$.  In particular, for $k=1$, $\mathcal{D}^{1,p}= G^{1,p}$, for all $p\in(1,\infty)$.
\end{thm}

\begin{defn}
Let $D^{*}$ denote the $L^2(\mu)$-adjoint of the derivative operator $D$, which has
domain in $L^{2}(\mathscr{W}\times[0,1],\mathscr{H})$ consisting of functions $G$ such that
\[
| \mathbb{E}[(DF,G)_{\mathscr{H}}]| \le C\|F\|_{L^{2}(\mu)},
\]
for all $F\in{\mathcal{D}}^{1,2}$, where $C$ is a constant depending on $G$.
For those functions $G$ in the domain of $D^{*}$, $D^{*}G$ is the element of
$L^{2}(\mu)$ such that
\[
\mathbb{E}[FD^{*}G] = \mathbb{E}[(DF,G)_{\mathscr{H}}].
\]
\end{defn}

It is known that $D$ is a continuous operator from $\mathcal{D}^{\infty}$ to
$\mathcal{D}^{\infty}(\mathscr{H})$, and similarly, $D^{\ast}$ is continuous from
$\mathcal{D}^{\infty}(\mathscr{H})$ to $\mathcal{D}^{\infty}$; see for example
%Proposition 1.5.4 from Nualart \cite{Nualart95}, or 
Theorem V-8.1 and its corollary in \cite{IkedaWatanabe89} .

\subsection{Lie group-valued Wiener functions}
\label{s:main}

Let $G$ be a Lie group with Lie algebra $\mathfrak{g}=\mathrm{Lie}(G)$ and identity element $e$, and let $L_g$ denote left translation by an element $g\in G$, and $R_g$ denote right 
translation.  Suppose $\{X_{i}\}_{i=1}^{k}\subset\mathfrak{g}$, and let $\mathfrak{g}_{0}:=\mathrm{span}(\{X_i\}_{i=1}^k)$.  Let $\left\langle  \cdot,\cdot \right\rangle$ be any inner product  on $\mathfrak{g}$, and extend $\left\langle  \cdot,\cdot\right\rangle$ to a right invariant metric on $G$ by defining $\left\langle \cdot,\cdot \right\rangle_g : T_gG \times T_gG \rightarrow \mathbb{R}$ as
\[
\left\langle v,w \right\rangle  _g := \left\langle  R_{g^{-1}*}v , R_{g^{-1}*}w\right\rangle  ,\quad\text{ for all }v,w\in T_gG.
\]
The $g$ subscript will be suppressed when there is no chance of confusion. 

\begin{notation}
Given an element $X\in\mathfrak{g}$, let $\tilde{X}$ denote the left invariant
vector field on $G$ such that $\tilde{X}(e)=X$, where $e$ is the identity of $G$.  Recall that $\tilde{X}$ left invariant means that the vector field commutes with left translation in the following way:
\[ \tilde{X}(f\circ L_g) = (\tilde{X} f) \circ L_g, \]
for all $f\in C^1(G)$.
Let $\hat{X}$ denote the right invariant vector field associated to $X$.
\end{notation}

\begin{notation}
Let $\operatorname{Ad}:G\rightarrow\mathrm{End}(\mathfrak{g})$ denote the adjoint representation of $G$ with differential $\operatorname{ad}:= d(\operatorname{Ad}):\mathfrak{g}\rightarrow\mathrm{End}(\mathfrak{g})$.  That is, $\operatorname{Ad}(g)=\operatorname{Ad}_g=L_{g*}R_{g^{-1}*}$, for all $g\in G$, and $\operatorname{ad}(X)=\operatorname{ad}_X=[X,\cdot]$, for all $X\in\mathfrak{g}$.  For any function $\varphi\in C^{1}(G)$, define $\hat{\nabla}\varphi,\tilde{\nabla}\varphi:G\rightarrow\mathfrak{g}$ such that, for any $g\in G$ and $X\in\mathfrak{g}$,
\[ \left\langle \hat{\nabla}\varphi(g),X \right\rangle := \left\langle d\varphi(g),R_{g*}X \right\rangle = (\hat{X}\varphi)(g) \]
and 
\[ \left\langle \tilde{\nabla}\varphi(g),X\right\rangle := \left\langle d\varphi(g),L_{g*}X\right\rangle  =(\tilde{X}\varphi)(g). \]
Then,
\[
\left\langle\tilde{\nabla}\hat{\nabla}\varphi(g) , X \otimes Y \right\rangle
   = \frac{d}{ds}\bigg|_{0} \frac{d}{dt}\bigg|_{0}\varphi \left(  e^{sX}ge^{tY} \right).
\]
for all $X,Y\in \mathfrak{g}$, and similarly for $\hat{\nabla}\tilde{\nabla}\varphi$.
\end{notation}

Now suppose $\{b_t^i\}_{i=1}^k$ are $k$ independent real-valued Brownian motions.  Then
\[ \vec{b}_t := X_i b_t^i := \sum_{i=1}^k X_i b^i_t \]
is a $(\mathfrak{g}_0,\langle\cdot,\cdot\rangle)$ Brownian motion.  In the sequel, the convention of summing over repeated upper and lower indices will be observed.  For  $h=(h^1,\ldots,h^k)\in \mathscr{H}$, let $\vec{h}:=X_ih^i\in\mathscr{H}(\mathfrak{g}_0)$.  Let $\xi:[0,1]\times \mathscr{W}\rightarrow G$ denote the solution to the Stratonovich stochastic differential equation
\begin{equation}
\label{e:SDE}
d\xi_{t} = \xi_t \circ d\vec{b}_t := L_{\xi_t*}\circ d\vec{b}_t 
   = L_{\xi_t*}X_i \circ db_t^i = \tilde{X}_i(\xi_t) \circ db_t^i,
\text{ with } \xi_0 = e.
\end{equation}
The solution $\xi$ exists by standard theory;  see, for example, Theorem V-1.1 of \cite{IkedaWatanabe89}. 

\begin{thm}
\label{t:xi.smooth}
For any $f\in C_c^\infty(G)$, $f(\xi_t)\in\mathcal{D}^\infty$ for all $t\in[0,1]$.  In particular, $D[f(\xi_t)]\in\mathscr{H}\otimes\mathbb{R}^k$ and 
\begin{equation*}
(D[f(\xi_t)])^i = \left\langle \hat{\nabla}f(\xi_t),  \int_0^{\cdot\wedge t}\operatorname{Ad}_{\xi_\tau} X_i \,d\tau\right\rangle,
\end{equation*}
for $i=1,\ldots,k$, componentwise in $\mathscr{H}$, and for any $h\in \mathscr{H}$,
\begin{equation*}
\partial_h [f(\xi_t)] = \int_0^t \left\langle \hat{\nabla}f(\xi_t), \operatorname{Ad}_{\xi_\tau} X_i \right\rangle \dot{h}^i_\tau \,d\tau.
\end{equation*}
\end{thm}

The proof of this theorem goes through a series of convergence arguments for solutions to two cutoff versions of Equation (\ref{e:SDE}).  It is necessary to first consider a version with vector fields modified by a cutoff function which is compactly supported on $G$ so that Taniguchi's result may be applied.  Another cutoff function which is compactly supported on $\mathrm{End}(\mathfrak{g})$ is then required to resolve certain convergence issues which arise with the derivatives (see Propositions \ref{p:d2.55} and \ref{p:d2.5} below).  Set the following notation for cutoff functions for the sequel.

\begin{notation}
\label{n:phim}
Let $|g|$ denote the distance from a point $g\in G$ to $e$ with respect to the right invariant metric.  Let $\{\varphi_m\}_{m=1}^{\infty}\subset C_{c}^{\infty}(G,[0,1])$ be a sequence of
functions, $\varphi_m\uparrow1$, such that $\varphi_m(g)=1$ when $\left|
g\right| \le m$ and $\sup_m\sup_{g\in G}|\tilde{\nabla}^{k}
\varphi_m(g)|<\infty$ for $k=0,1,2,\ldots$; see Lemma 3.6 of Driver and Gross
\cite{Driver92b}. 

Fix $\psi\in C_c^\infty(\mathrm{End}(\mathfrak{g}),[0,1])$ so that $\psi=1$ near $I$ and $\psi(x)=0$ if $|x|\ge2$, where $|\cdot|$ is the distance from $I$ with respect to any metric on $\mathrm{End}(\mathfrak{g})$, and define
\[ \left\langle\psi'(x),A\right\rangle := \frac{d}{d\varepsilon}\bigg|_{\varepsilon=0} \psi(x+\varepsilon A), \]  
for any $A\in\mathrm{End}(\mathfrak{g})$.  Take $v(g) := \psi(\operatorname{Ad}_g)$ and $u_m(g) := \varphi_m(g) v(g)$.  
\end{notation}

\section{Proof of Theorem \ref{t:xi.smooth}}
\label{s:SDEp}
There are several standard convergence results on matrix groups which are assumed in the sequel.  For the sake of completeness (and since they do not seem to exist in one place in the current literature), these are compiled in the note \cite{app}.

The following proposition will be used repeatedly; see for example Driver \cite{Driver04}. 
\begin{prop}
\label{p:d4.1}
Suppose $p\in[2,\infty)$, $\alpha_t$ is a predictable $\mathbb{R}^d$--valued process,
$A_t$ is a predictable $\mathrm{Hom}(\mathfrak{g}_0,\mathbb{R}^d)$--valued process, and
\begin{equation*}
Y_t := \int_0^t A_{\tau} \,d\vec{b}_{\tau} + \int_0^t \alpha_{\tau} \,d\tau = \int_0^t A_{\tau} X_i\,db^i_{\tau} + \int_0^t \alpha_{\tau} \,d\tau,
\end{equation*}
where $\{b^1,\ldots,b^k\}$ are $k$ independent real Brownian motions.  Then
\begin{equation*}
\mathbb{E} \sup_{\tau\le t}\left| Y_{\tau}\right|^p \le C_{p} \left\{  \mathbb{E} \left(  \int_0^t\left| A_{\tau}\right| ^2d\tau\right)^{p/2}
    + \mathbb{E} \left(  \int_0^t\left| \alpha_{\tau}\right| d\tau\right)^p\right\},
\end{equation*}
where
\[
|A|^2 = \operatorname{tr}(AA^*) = \sum_{i=1}^n (AA^*)_{ii} = \sum_{i,j}A_{ij}A_{ij} = \operatorname{tr}(A^*A).
\]
\end{prop}

\begin{notation}
Let $\delta_n$ denote constants such that $\lim_{n\rightarrow\infty}\delta_n=0$.  Also, write $f\lesssim g$, if there is a positive constant $C$ so that $f\le Cg$.
\end{notation}

\begin{lem}
\label{l.d2.1}
Let $u\in C^1(G)$ such that $u$ and $\tilde{X} u$ are bounded for all $X\in\mathfrak{g}_0$.
Then the solution $\eta:[0,1]\times \mathscr{W}\rightarrow G$ to the stochastic differential equation
\[
d\eta = u(\eta)\eta \circ d\vec{b} := u(\eta)L_{\eta*} \circ d\vec{b} = u(\eta) \tilde{X}_i (\eta) \circ db^i,
\text{ with }\eta_{0}=e
\]
exists for all time; that is, $\eta$ has no explosion.
\end{lem}

\begin{proof}
Let $\zeta$ be the life-time of $\eta$ and $\varphi\in C_c^\infty(G)$. Then on $\{t<\zeta\}$,
\begin{align}
d\varphi(&\eta) =  u(\eta) \tilde{X}_i \varphi(\eta) \circ db^i
   = u(\eta)  \left\langle \tilde{\nabla}\varphi(\eta)  ,X_i\right\rangle  \circ db^i
   = u\left(\eta\right)  \left\langle  \tilde{\nabla}\varphi(\eta) , \circ d\vec{b}\right\rangle \notag \\
   &= u(\eta)  \left\langle  \tilde{\nabla}\varphi(\eta),d\vec{b}\right\rangle  
     + \frac{1}{2}d\left[  u(\eta)  \left\langle \tilde{\nabla}\varphi(\eta)  ,\cdot\right\rangle  \right]  d\vec{b} \notag \\
   &= u(\eta)  \left\langle  \tilde{\nabla}\varphi(\eta),d\vec{b}\right\rangle  
     + \frac{1}{2} \left\langle \tilde{\nabla} \left[ u(\eta) \left\langle \tilde{\nabla}\varphi(\eta)  ,\cdot\right\rangle \right]  ,d\vec{b}\right\rangle  d\vec{b} \notag \\
   &=u(\eta) \left\langle \tilde{\nabla} \varphi(\eta),d\vec{b} \right\rangle + \frac{1}{2}\left[\left\langle  \tilde{\nabla}u(\eta),d\vec{b} \right\rangle \left\langle \tilde{\nabla} \varphi(\eta),d\vec{b} \right\rangle  
     + u(\eta) \left\langle \tilde{\nabla}^2 \varphi(\eta)  ,d\vec{b}\otimes d\vec{b}\right\rangle \right] \notag \\
 \label{e:d2.2}
   &=u(\eta)  \left\langle  \tilde{\nabla}\varphi(\eta),d\vec{b}\right\rangle  \notag\\
   & \qquad +\frac{1}{2}\sum_{i=1}^k\left[  \left\langle  \tilde{\nabla} u(\eta),X_i \right\rangle \left\langle  \tilde{\nabla} \varphi(\eta),X_i\right\rangle  
      + u(\eta)  \left\langle  \tilde{\nabla}^2\varphi(\eta) , X_i\otimes X_i \right\rangle  \right]  dt.
\end{align}

Let $\{\varphi_m\}_{m=1}^{\infty}\subset C_{c}^{\infty}(G,[0,1])$ be a sequence of
functions as in Notation \ref{n:phim}.  Then from Equation (\ref{e:d2.2}) (using the convention that
$\varphi_m(\eta_t) = 0$, $\tilde{\nabla}\varphi_m(\eta_t) = 0$, and $\tilde{\nabla}^2\varphi_m(\eta_t) = 0$ 
on $\{t>\zeta\}$)
\begin{multline*}
\varphi_m(\eta_t) = 1 + \int_0^{t} u(\eta) \left\langle \tilde{\nabla}\varphi_m(\eta)  ,d\vec{b}\right\rangle  \\
   + \frac{1}{2}\int_0^{t}\sum_{i=1}^k\left[  \left\langle  \tilde{\nabla}u(\eta) , X_i \right\rangle 
   \left\langle \tilde{\nabla}\varphi_m(\eta) , X_i \right\rangle + u(\eta) \left\langle \tilde{\nabla}^2 \varphi_m(\eta) , X_i\otimes X_i \right\rangle\right]  d\tau.
\end{multline*}
Taking the expectation of this equation then gives
\[
\mathbb{E} \left[\varphi_m(\eta_t)\right] = 1 + \delta_m,
\]
where
\[
\delta_m := \frac{1}{2}\mathbb{E} \int_0^{t}  \sum_{i=1}^k \left[  \left\langle \tilde{\nabla}u(\eta) , X_i \right\rangle 
   \left\langle \tilde{\nabla}\varphi_m(\eta) , X_i \right\rangle + u(\eta) \left\langle \tilde{\nabla}^2 \varphi_m(\eta) , X_i\otimes X_i \right\rangle \right]  d\tau.
\]
Now by construction of the $\varphi_m$ and the assumptions on $u$, there is a constant $M<\infty$ such that
\[
\left|  \left\langle  \tilde{\nabla}u(\eta_\tau) ,X_i \right\rangle  
   \left\langle  \tilde{\nabla}\varphi_m(\eta_\tau) , X_i \right\rangle + u(\eta_\tau) \left\langle \tilde{\nabla}^2 \varphi_m(\eta_\tau) , 
   X_i \otimes X_i \right\rangle \right| \le M.
\]
Moreover, $\lim_{m\rightarrow\infty} \left| \tilde{\nabla}\varphi_m \right| = 0 = \lim_{m\rightarrow\infty}\left| \tilde{\nabla}^2\varphi_m\right|$.
So it follows by the dominated convergence theorem that $\lim_{m\rightarrow\infty} \delta_m = 0$, and 
\[
1 = \lim_{m\rightarrow\infty}\mathbb{E}\left[ \varphi_m(\eta_t)\right] = \mathbb{E} \left[\lim_{m\rightarrow\infty} \varphi_m(\eta_t) \right]
   = \mathbb{E} 1_{\{t<\zeta\}} = P(t<\zeta);
\]
that is, $\zeta=\infty$ $\mu$-a.s.
\end{proof} 

%% ALTERNATELY,
\iffalse
\begin{proof}
Let $|g|$ denote the distance from a point $g\in G$ to $e$ with respect to the right invariant metric.
Let $\{\varphi_m\}_{m=1}^{\infty}\subset C_{c}^{\infty}(G,[0,1])$ be a sequence of
functions, $\varphi_m\uparrow1$, such that $\varphi_m(g)=1$ when $\left|
g\right| \le m$ and $\sup_m\sup_{g\in G}|\tilde{\nabla}^{k}
\varphi_m(g)|<\infty$ for $k=0,1,2,\ldots$; see Lemma 3.6 of Driver and Gross
\cite{Driver92b}.  Also, let $\{K_n\}_{n=1}^\infty$ be a sequence of nested compact sets in $G$, such that $K_n\uparrow G$, and take $\tau_n$ to be the exit time from $K_n$,
\[ \tau_n := \inf\{t>0: \eta_t\notin K_n \}. \]
Let $\eta_t^n:=\eta_{t\wedge\tau_n}$.  Then, by construction of $\varphi_m$ and $\tau_n$,
\[ 
1 = \lim_{m\rightarrow\infty}\mathbb{E}\left[ \varphi_m(\eta^n_t)\right] = \mathbb{E} \left[\lim_{m\rightarrow\infty} \varphi_m(\eta^n_t) \right]
   = \mathbb{E} 1_{\{t\wedge\tau_n<\zeta\}} = P(t\wedge\tau_n<\zeta),
\]
for all $n=1,2,\ldots.$  Thus, $1=\lim_{n\rightarrow\infty} P(t\wedge\tau_n<\zeta) = P(t<\zeta)$, and so $\zeta=\infty$ $\mu$-a.e.
\end{proof}
\fi
%% must switch the order of limits here
%%

Now let $u\in C_c^\infty(G)$, and suppose $\eta:[0,1]\times \mathscr{W}\rightarrow G$ is a solution to the stochastic differential equation
\begin{equation}
\label{e:gSDE}
d\eta = u(\eta)\eta\circ d\vec{b} := u(\eta)L_{\eta*}\circ d\vec{b} =  u(\eta)\tilde{X}_i(\eta) \circ db^{i}.
\end{equation}
The previous lemma implies that the solution $\eta$ exists for all time.  Since the vector fields $u\tilde{X}_i $ have compact support, $G$ may be embedded as a Euclidean submanifold in a ``nice" way so that the embedded vector fields are bounded with bounded derivatives.  Then Theorem 2.1 of Taniguchi \cite{Taniguchi83} implies that, for any $f\in C_c^\infty(G)$, the function $f(\eta_t)\in\mathcal{D}^\infty$.

\begin{prop}
\label{p:Dh}
Fix $h\in \mathscr{H}$, and let $\eta$ be the solution to Equation (\ref{e:gSDE}).  Then, for $f\in C_c^\infty(G)$,
\[ \partial_h f(\eta_t) = \left\langle \hat{\nabla}f(\eta_t), \theta_t \right\rangle, \]
where $\theta_t:\mathscr{W}\rightarrow\mathfrak{g}$ solves
\begin{equation}
\label{e:d2.z}
\theta_t = \int_0^t \left(\left\langle \hat{\nabla}u(\eta),\theta \right\rangle \operatorname{Ad}_{\eta} \circ d\vec{b} + u(\eta)\operatorname{Ad}_{\eta} d\vec{h} \right).
\end{equation}
Writing this in It\^o form gives
\begin{multline*}
\theta_t = \int_0^t \bigg(\left\langle\hat{\nabla}u(\eta) , \theta\right\rangle \operatorname{Ad}_\eta d\vec{b} + u(\eta) \operatorname{Ad}_\eta d\vec{h} \\
   + \frac{1}{2}\sum_{i=1}^k\left\langle\tilde{\nabla}\hat{\nabla}u(\eta) , X_i \otimes\theta\right\rangle \operatorname{Ad}_\eta X_i \,d\tau\bigg).
\end{multline*}
\end{prop}

\begin{proof}
Let $\eta_t^s :=  \eta_t( \cdot+sh):[0,1]^2\times \mathscr{W}\rightarrow G$.  Then $\eta_t^s$ satisfies the Stratonovich equation in $t$,
\begin{equation}
\label{e:pseudo}
d\eta_t^s = u(\eta_t^s)L_{\eta_t^s*} \left(\circ d\vec{b}_t + sd\vec{h}_t\right).
\end{equation}
By Corollary 4.3 of Driver \cite{Driver92a}, there exists a modification of $\eta_t^s$ so that the mapping $s\mapsto\eta_t^s$ is smooth in the sense that, for any function $f\in C^\infty_c(G)$, $s\mapsto f(\eta_t^s)$ is smooth, and, furthermore, for any one-form $\vartheta$ acting on $T_{\eta_t^s}G$,
\begin{equation}
\label{e:qq} 
\frac{\partial}{\partial s}\bigg|_0 \int_0^t \vartheta(d\eta_\tau^s) = \int_0^t d\vartheta \left(\frac{\partial}{\partial s}\bigg|_0\eta_\tau^s, 
   d\eta_\tau\right) + \vartheta\left(\frac{\partial}{\partial s}\bigg|_0 \eta_\tau^s\right) \bigg|_{\tau=0}^t. 
\end{equation}

Let $\vartheta$ be the $\mathfrak{g}$-valued one-form such that $\vartheta (\hat{X})=X$, for all $X\in\mathfrak{g}$.  Since $\vartheta$ is right invariant, the two-form $d\vartheta$
%: T_{ g^s_t}G\times T_{ g^s_t}G \rightarrow \mathfrak{g}$
satisfies the identity
\begin{equation*}
d\vartheta = \vartheta\wedge\vartheta
\end{equation*}
where $\vartheta\wedge\vartheta (X,Y):=[\vartheta(X),\vartheta(Y)]$ for any $X,Y\in\mathfrak{g}$; see for example \cite{AgricolaFriedrich}.  
%page 210
Let $\theta_t:= \vartheta\left(\frac{\partial}{\partial s}\big|_0\eta_t^s\right)$, so that $\frac{\partial}{\partial s}\big|_0 \eta_t^s = R_{\eta_t*} \theta_t$.  Thus,
\[ \frac{\partial}{\partial s}\bigg|_0 \operatorname{Ad}_{\eta_t^s} = d\operatorname{Ad}(\theta_t) \operatorname{Ad}_{\eta_t} = \operatorname{ad}_{\theta_t} \operatorname{Ad}_{\eta_t}. \]
By Equation (\ref{e:qq}), 
\begin{align*} 
\theta_t &= \vartheta\left(\frac{\partial}{\partial s}\bigg|_0 \eta_t^s\right) = \frac{\partial}{\partial s}\bigg|_0 \int_0^t \vartheta(d\eta_\tau^s)
      - \int_0^t \left[\vartheta\left(\frac{\partial}{\partial s}\bigg|_0 \eta_t^s\right)  , \vartheta\left(d\eta_\tau\right)\right] \\
   &= \frac{\partial}{\partial s}\bigg|_0 \int_0^t u(\eta_t^s)\operatorname{Ad}_{\eta_t^s} \left(\circ d\vec{b}_\tau + sd\vec{h}_\tau \right) 
      - \int_0^t \left[\theta_\tau  , u(\eta_t)\operatorname{Ad}_{\eta_t} \circ d\vec{b}_\tau \right] \\
   &= \int_0^t \left(\left\langle \hat{\nabla}u(\eta_\tau),\theta_\tau \right\rangle \operatorname{Ad}_{\eta_\tau} \circ 
      d\vec{b}_\tau + u(\eta_\tau)\operatorname{ad}_{\theta_\tau} \operatorname{Ad}_{\eta_\tau} \circ d\vec{b}_\tau  
      + u(\eta_t)\operatorname{Ad}_{\eta_t} d\vec{h}_\tau \right) \\
   & \qquad   - \int_0^t \left[\theta_\tau , u(\eta_t)\operatorname{Ad}_{\eta_t} \circ d\vec{b}_\tau\right] \\
   &=  \int_0^t \left(\left\langle \hat{\nabla}u(\eta_\tau),\theta_\tau \right\rangle \operatorname{Ad}_{\eta_\tau} \circ d\vec{b}_\tau  
      + u(\eta_\tau)\operatorname{Ad}_{\eta_\tau} d\vec{h}_\tau \right),
\end{align*}
and for any $f\in C_c^\infty(G)$,
 \[ \partial_h f(\eta_t) = \frac{\partial}{\partial s} \bigg|_0 f(\eta_t^s) = \left\langle \hat{\nabla}f(\eta_t), \vartheta\left(\frac{\partial}{\partial s}\bigg|_0 \eta_t^s\right) \right\rangle
    =\left\langle \hat{\nabla} f(\eta_t), \theta_t \right\rangle. \]

Now, to write this equation in It\^o form, first note that 
\begin{align}
d\operatorname{Ad}_\eta  
   &= u(\eta) \operatorname{Ad}_\eta \circ\operatorname{ad}_{d\vec{b}} \notag \\
   &= u(\eta) \operatorname{Ad}_\eta \operatorname{ad}_{d\vec{b}} + \frac{1}{2}\left[ \left\langle \tilde{\nabla}u(\eta) , d\vec{b} \right\rangle \operatorname{Ad}_\eta
     + u^2(\eta)  \operatorname{Ad}_\eta \operatorname{ad}_{d\vec{b}}\right]  \cdot \operatorname{ad}_{d\vec{b}} \notag \\
\label{e:d2.4}
   &= u(\eta)  \operatorname{Ad}_\eta \operatorname{ad}_{d\vec{b}}+\frac{1}{2}\sum_{i=1}^k\left[ \left\langle\tilde{\nabla}u(\eta)  ,X_i\right\rangle \operatorname{Ad}_\eta \operatorname{ad}_{X_i}
     +u^2(\eta) \operatorname{Ad}_\eta \operatorname{ad}_{X_i}^2 \right] dt,
\end{align}
where $\operatorname{ad}_{d\vec{b}}=\operatorname{ad}_{X_i} db^i$.  This then implies
\begin{align*}
d\bigg[  \left\langle\hat{\nabla}u(\eta)  ,\theta\right\rangle &\operatorname{Ad}_\eta\bigg]\cdot d\vec{b}  
   =\left[  \left\langle\tilde{\nabla}\hat{\nabla}u(\eta),d\vec{b}\otimes\theta\right\rangle \operatorname{Ad}_\eta
     +  \left\langle\hat{\nabla}u(\eta),\theta\right\rangle u(\eta)  \operatorname{Ad}_\eta\operatorname{ad}_{d\vec{b}}\right]  d\vec{b}  \\
   &=\sum_{i=1}^k\left[  \left\langle\tilde{\nabla}\hat{\nabla}u(\eta) , X_i\otimes\theta\right\rangle \operatorname{Ad}_\eta
     +\left\langle \hat{\nabla}u(\eta) , \theta \right\rangle u(\eta) \operatorname{Ad}_\eta\operatorname{ad}_{X_i}\right]  X_i \,dt \\
   &=\sum_{i=1}^k\left\langle\tilde{\nabla}\hat{\nabla}u(\eta),X_i\otimes\theta\right\rangle \operatorname{Ad}_\eta X_i \,dt,
\end{align*}
using that $\operatorname{ad}_{X}X = [X,X] = 0$.  Therefore,
\begin{align*}
\left\langle\hat{\nabla}u(\eta)  ,\theta\right\rangle \operatorname{Ad}_\eta\circ d\vec{b}  
   &= \left\langle \hat{\nabla}u(\eta)  ,\theta\right\rangle \operatorname{Ad}_\eta d\vec{b} + \frac{1}{2} d\left[  \left\langle\hat{\nabla}u(\eta)  ,\theta\right\rangle
     \operatorname{Ad}_\eta \right]  \cdot d\vec{b} \\
   &= \left\langle \hat{\nabla}u(\eta)  ,\theta\right\rangle \operatorname{Ad}_\eta d\vec{b} + \frac{1}{2} \sum_{i=1}^k\left\langle\tilde{\nabla}\hat{\nabla}u(\eta)
     ,X_i \otimes \theta \right\rangle \operatorname{Ad}_\eta X_i \,dt,
\end{align*}
and $\theta$ satisfies the It\^o differential equation
\begin{equation*}
d\theta = \left\langle\hat{\nabla}u(\eta) , \theta\right\rangle \operatorname{Ad}_\eta d\vec{b} + u(\eta) \operatorname{Ad}_\eta d\vec{h} 
   + \frac{1}{2}\sum_{i=1}^k\left\langle\tilde{\nabla}\hat{\nabla}u(\eta) , X_i \otimes\theta\right\rangle \operatorname{Ad}_\eta X_i \,dt. 
\end{equation*}
\end{proof}

\begin{remark}
For $u_m=\varphi_m v$, the following derivative formulae hold:
\[
\hat{\nabla}u_m  
   =\hat{\nabla}v\cdot\varphi_m + v\cdot\hat{\nabla}\varphi_m \]
and
\begin{align*}
\tilde{\nabla}\hat{\nabla} u_m  
   &=\tilde{\nabla} \left[  \hat{\nabla} v \cdot \varphi_m + v \cdot \hat{\nabla}\varphi_m \right] \\
   &=\tilde{\nabla}\hat{\nabla}v\cdot\varphi_m+\tilde{\nabla}v\otimes
        \hat{\nabla}\varphi_m+\tilde{\nabla}\varphi_m\otimes\hat{\nabla}v
        + v \cdot\tilde{\nabla}\hat{\nabla}\varphi_m. 
 \end{align*}
\end{remark}

\begin{prop}
\label{p:eta}
Let $\eta^m:[0,1]\times \mathscr{W}\rightarrow G$ denote the solution to the equation 
\begin{equation}
\label{e:d2.7}
\begin{split}
d\eta^m = u_m(\eta^m) \eta^m \circ d\vec{b} 
   &= \varphi_m(\eta^m)v(\eta^m) \eta^m \circ d\vec{b} \\
   &= \varphi_m(\eta^m)\psi(\operatorname{Ad}_{\eta^m}) \eta^m \circ d\vec{b}, \text{ with } \eta_0^m = e,
\end{split}
\end{equation}
and $\eta:[0,1]\times \mathscr{W}\rightarrow G$ denote the solution to
\begin{equation}
\label{e:d2.8}
d\eta = v(\eta)\eta \circ d\vec{b} = \psi(\operatorname{Ad}_\eta) \eta \circ d\vec{b}, \text{ with } \eta_{0} = e.
\end{equation}
Then for all $f\in C_c^\infty(G)$,
\[ \lim_{m\rightarrow\infty} \mathbb{E}\sup_{\tau\le1} |f(\eta_\tau^m) - f(\eta_\tau)|^p = 0, \]
for all $p\in(1,\infty)$.
\end{prop}

\begin{proof}
By Lemma \ref{l.d2.1}, Equation (\ref{e:d2.7}) has a global solution.  Notice also that $\tilde{X} v(\eta) = \left\langle \psi'(\operatorname{Ad}_\eta),\operatorname{ad}_X\right\rangle$ is bounded, so that Equation(\ref{e:d2.8}) also has a global solution.  Then by Equation (\ref{e:d2.2}),
\begin{multline*}
df(\eta^m) = u_m(\eta^m)  \left\langle  \tilde{\nabla}f(\eta^m),d\vec{b}\right\rangle  \\
   +\frac{1}{2}\sum_{i=1}^k\left[  \left\langle  \tilde{\nabla} u_m(\eta^m),X_i \right\rangle \left\langle  \tilde{\nabla} f(\eta^m),X_i\right\rangle  
   + u_m(\eta^m)  \left\langle  \tilde{\nabla}^2f(\eta^m) , X_i\otimes X_i \right\rangle  \right]  dt,
\end{multline*}
and similarly,
\begin{multline*}
df(\eta) =v(\eta)  \left\langle  \tilde{\nabla}f(\eta),d\vec{b}\right\rangle  \\
   + \frac{1}{2}\sum_{i=1}^k\left[  \left\langle  \tilde{\nabla} v(\eta),X_i \right\rangle \left\langle  \tilde{\nabla} f(\eta),X_i\right\rangle  
   + v(\eta)  \left\langle  \tilde{\nabla}^2f(\eta) , X_i\otimes X_i \right\rangle  \right]  dt.
\end{multline*}
Thus,
\begin{multline*} 
d[f(\eta^m)-f(\eta)] = u_m(\eta^m)  \left\langle  \tilde{\nabla}f(\eta^m),d\vec{b}\right\rangle - v(\eta)\left\langle\tilde{\nabla}f(\eta),d\vec{b}\right\rangle  \\
    + \frac{1}{2}\sum_{i=1}^k\left[  \left\langle  \tilde{\nabla} u_m(\eta^m),X_i \right\rangle \left\langle  \tilde{\nabla} f(\eta^m),X_i\right\rangle 
      - \left\langle  \tilde{\nabla} v(\eta),X_i \right\rangle \left\langle  \tilde{\nabla} f(\eta),X_i\right\rangle \right. \\
    \left. + u_m(\eta^m)  \left\langle  \tilde{\nabla}^2f(\eta^m) , X_i\otimes X_i \right\rangle 
      - v(\eta)  \left\langle  \tilde{\nabla}^2f(\eta) , X_i\otimes X_i \right\rangle \right]  dt.
\end{multline*}   
Bound this expression by applying Proposition \ref{p:d4.1} to each term.  For the first term, note that
$u_m\rightarrow v$ boundedly, as $m\rightarrow\infty$, and this implies that
\begin{align*}
\mathbb{E} \bigg| \int_0^t  \bigg[ u_m(\eta^m)  \big\langle \tilde{\nabla}f(&\eta^m),d\vec{b}\big\rangle 
      - v(\eta)\big\langle\tilde{\nabla}f(\eta),d\vec{b}\big\rangle \bigg] \bigg|^p \\
   &\lesssim \sum_{i=1}^k \mathbb{E} \int_0^t \left|u_m(\eta^m)  \left\langle  \tilde{\nabla}f(\eta^m),X_i\right\rangle 
      - v(\eta)\left\langle\tilde{\nabla}f(\eta),X_i\right\rangle \right|^p d\tau \\ 
   &\lesssim \|\tilde{\nabla}f\|_\infty \mathbb{E} \int_0^t \left|u_m(\eta^m) - v(\eta)\right|^p d\tau \rightarrow 0,
\end{align*}
as $m\rightarrow\infty$, by the dominated convergence theorem.  Similarly, for the second term, $\tilde{\nabla}u_m\rightarrow\tilde{\nabla}v$ boundedly, as $m\rightarrow\infty$, implies that
\begin{align*}
 \mathbb{E} \bigg| \int_0^t \bigg[  &\left\langle  \tilde{\nabla} u_m(\eta^m),X_i \right\rangle \left\langle  \tilde{\nabla} f(\eta^m),X_i\right\rangle 
      - \left\langle  \tilde{\nabla} v(\eta),X_i \right\rangle \left\langle  \tilde{\nabla} f(\eta),X_i\right\rangle \bigg] d\tau \bigg|^p \\
   &\lesssim \mathbb{E}\int_0^t  \left|\left\langle  \tilde{\nabla} u_m(\eta^m),X_i \right\rangle \left\langle  \tilde{\nabla} f(\eta^m),X_i\right\rangle 
      - \left\langle  \tilde{\nabla} v(\eta),X_i \right\rangle \left\langle  \tilde{\nabla} f(\eta),X_i\right\rangle \right|^p d\tau \\
   &\lesssim \|\tilde{\nabla}f\|_\infty \mathbb{E} \int_0^t \left|\tilde{\nabla}u_m(\eta^m) 
      - \tilde{\nabla}v(\eta)\right|^p d\tau
    \rightarrow 0,
\end{align*}
as $m\rightarrow\infty$, by the dominated convergence theorem.  Finally,
\begin{align*}
 \mathbb{E} \bigg|\int_0^t \bigg[ u_m(&\eta^m)  \left\langle  \tilde{\nabla}^2f(\eta^m) , X_i\otimes X_i \right\rangle 
      - v(\eta)  \left\langle  \tilde{\nabla}^2f(\eta) , X_i\otimes X_i \right\rangle \bigg] d\tau \bigg|^p \\
   &= \mathbb{E} \int_0^t \left| u_m(\eta^m) \left\langle  \tilde{\nabla}^2f(\eta^m) , X_i\otimes X_i \right\rangle 
      - v(\eta)  \left\langle  \tilde{\nabla}^2f(\eta) , X_i\otimes X_i \right\rangle \right|^p d\tau \\
   &\lesssim \|\tilde{\nabla}^2f\|_\infty  \mathbb{E} \int_0^t \left|u_m(\eta^m) - v(\eta)\right|^p d\tau
   \rightarrow 0,
\end{align*}
as $m\rightarrow\infty$, again by dominated convergence.  Thus,
\[ \lim_{m\rightarrow\infty} \mathbb{E}\sup_{\tau\le1} |f(\eta_\tau^m)-f(\eta_\tau)|^p = 0, \]
as desired.
\end{proof}

\begin{prop}
\label{p:d2.55}
Let $U^m_t = \operatorname{Ad}_{\eta^m_t}:\mathscr{W}\rightarrow \mathrm{End}(\mathfrak{g})$ and $U_t = \operatorname{Ad}_{\eta_t}:\mathscr{W}\rightarrow\mathrm{End}(\mathfrak{g})$, which satisfy the stochastic differential equations
\begin{equation}
\label{e:d2.11}
dU^m = u_m(\eta^m) U^m \circ\operatorname{ad}_{d\vec{b}} = \varphi_m(\eta^m)\psi(U^m) U^m \circ\operatorname{ad}_{d\vec{b}}, 
   \text{ with } U_0^m = I,
\end{equation}
and
\begin{equation}
\label{e:d1}
dU = v(\eta) U \circ\operatorname{ad}_{d\vec{b}} = \psi(U) U \circ\operatorname{ad}_{d\vec{b}}, \text{ with } U_0 = I,
\end{equation}
where $\operatorname{ad}_{d\vec{b}}=\operatorname{ad}_{X_i} db^i$.  Then
\[
\lim_{m\rightarrow\infty} \mathbb{E} \sup_{\tau\le1} | U^m_\tau - U_\tau |^p = 0,
\]
for all $p\in(1,\infty)$.
\end{prop}

\iffalse
\begin{proof}
For any representation $\rho:G\rightarrow \Aut(V)$ of $G$, $g\in G$, and $X\in\mathfrak{g}$,
\[ \begin{split}
\rho_*(R_{g*}X) = \de\bigg|_0 \rho(g\cdot e^{\e X}) =  \de\bigg|_0 \rho(g)\rho( e^{\e X}) = \rho(g)d\rho(X) \\
\rho_*(L_{g*}X) = \de\bigg|_0 \rho(e^{\e X}\cdot g) =  \de\bigg|_0\rho( e^{\e X}) \rho(g) = d\rho(X) \rho(g),
\end{split} \]
where we have used that $\rho$ is a homomorphism  in the second equality.  Then 
\[ \begin{split} 
\frac{d}{dt}\rho( g_t) &= \varphi(g_t)d\rho(d\vec{b}_t) \rho( g_t) \\
\rho( g_0) &= I
\end{split} \]
where $d\rho(Y)= \frac{d}{dt}\big|_0 \rho(e^{tY})\in\mathrm{End}(\mathfrak{g})$ for any $Y\in\mathfrak{g}$.
Thus, for $\rho=\operatorname{Ad}:G\rightarrow\Aut(\mathfrak{g})$, $d\rho=\operatorname{ad}:\mathfrak{g}\rightarrow\mathrm{End}(\mathfrak{g})$,
and $\operatorname{Ad}_{ g_t}:\mathscr{W}\rightarrow \Aut(\mathfrak{g})$ satisfies the following stochastic differential
equation:
\begin{equation}
\begin{split}
d\operatorname{Ad}_{ g_t} &=  \varphi(g_t)\circ \operatorname{ad}_{d\vec{b}_t} \operatorname{Ad}_{ g_t} =  \varphi(g_t)\operatorname{ad}_{X_i}\operatorname{Ad}_{ g_t} \circ db_t^i \\
\operatorname{Ad}_{g_0} &= I,
\end{split}
\end{equation}
where $\circ\operatorname{ad}_{d\vec{b}_t}=\sum_{i=1}^k \operatorname{ad}_{X_i}\circ db^i_t$.
%\[ \begin{split}
%d\rho( g_t) = ``d\rho dg_t '' = d\rho L_{g_t*} db_t =
%   &= \de\bigg|_0 \rho( g_t e^{\e\circ d b_t}) \\
%   &= \de\bigg|_0 \rho( g_t) \rho(e^{\e\circ d b_t})
%   = \rho( g_t) \circ d\rho(d b_t),
%\end{split} \]
\end{proof}
\fi

\begin{proof}
Using Equation (\ref{e:d2.4}), rewrite(\ref{e:d2.11}) and (\ref{e:d1}) in It\^o form as
\begin{align*}
d&U^m 
   = u_m(\eta^m)U^m\operatorname{ad}_{d\vec{b}} + \frac{1}{2}\sum_{i=1}^k\left[ \left\langle\tilde{\nabla}u_m(\eta^m)  ,
     X_i\right\rangle U^m \operatorname{ad}_{X_i} + u_m^2(\eta^m) U^m \operatorname{ad}_{X_i}^2 \right] dt \\
   &=  \varphi_m(\eta^m)\psi(U^m) U^m \operatorname{ad}_{d\vec{b}} + \frac{1}{2}\sum_{i=1}^k
      \bigg[ \left\langle \tilde{\nabla}\varphi_m(\eta^m), {X_i}\right\rangle \psi(U^m) U^m \operatorname{ad}_{X_i} \\
   & \qquad + \varphi_m(\eta^m)\psi(U^m) \left\langle \psi'(U^m) , U^m\operatorname{ad}_{X_i} \right\rangle U^m \operatorname{ad}_{X_i} 
     + \varphi_m^2(\eta^m)\psi^2(U^m) U^m \operatorname{ad}_{X_i}^2 \bigg] dt
\end{align*}
and
\begin{align*}
dU 
%   &=  \psi(U) U \operatorname{ad}_{db} + \frac{1}{2}d\left[\psi(U) U \right]  \cdot \operatorname{ad}_{db} \\
%   &=  \psi(U) U \operatorname{ad}_{db} + \frac{1}{2}\left[ \psi(U))\left\langle\psi'(U),U\operatorname{ad}_{db}\right\rangle U + \psi^2(U)U\operatorname{ad}_{db} \right] \cdot 
%   \operatorname{ad}_{db} \\
   &= v(\eta)  U \operatorname{ad}_{d\vec{b}} + \frac{1}{2}\sum_{i=1}^k\left[ \left\langle\tilde{\nabla}v(\eta)  ,X_i\right\rangle U \operatorname{ad}_{X_i}
     +v^2(\eta) U \operatorname{ad}_{X_i}^2 \right] dt \\
   &= \psi(U) U \operatorname{ad}_{d\vec{b}} + \frac{1}{2} \sum_{i=1}^k
      \left[ \psi(U) \left\langle \psi'(U), U\operatorname{ad}_{X_i} \right\rangle U\operatorname{ad}_{X_i} + \psi^2(U) U \operatorname{ad}^2_{X_i} \right] dt.
\end{align*}
Let $\varepsilon^m := U^m - U$.  Then by the above,
\begin{multline}
\label{e:zs}
d\varepsilon^m 
   = \left( \varphi_m(\eta^m)\psi(U^m) U^m - \psi(U) U \right) \operatorname{ad}_{d\vec{b}} \\ 
   + \frac{1}{2} \sum_{i=1}^k \bigg[ \left\langle \tilde{\nabla}\varphi_m(\eta^m), {X_i}\right\rangle 
      \psi(U^m) U^m \operatorname{ad}_{X_i} \\
   + \left(\varphi_m(\eta^m)\psi(U^m) \left\langle \psi'(U^m) , U^m\operatorname{ad}_{X_i} \right\rangle U^m 
   - \psi(U)\left\langle \psi'(U), U\operatorname{ad}_{X_i} \right\rangle U \right)\operatorname{ad}_{X_i} \\
   + \left( \varphi_m^2(\eta^m)\psi^2(U^m) U^m -  \psi^2(U) U \right) \operatorname{ad}^2_{X_i} \bigg] dt.
\end{multline}
Again apply Proposition \ref{p:d4.1} to work term by term to bound the above expression.  Note first that, since $\psi$ has compact support, $U$ and $U^m$ always remain in a {\em fixed} compact
subset of $\mathrm{End}(\mathfrak{g})$.  Thus,
\begin{align*}
\mathbb{E} \bigg| \int_0^t  [ \varphi_m(&\eta^m) \psi(U^m)U^m - \psi(U)U ]  \operatorname{ad}_{d\vec{b}}\bigg| ^p \\
   &\lesssim \mathbb{E} \int_0^t \left| \varphi_m(\eta^m)\psi(U^m)U^m - \psi(U)U \right| ^p d\tau \\
   &\lesssim \mathbb{E} \int_0^t \left| \psi(U^m) U^m - \psi(U)U  \right|^p d\tau 
      + \mathbb{E} \int_0^t |\varphi_m(\eta^m) - 1|^p |\psi(U^m)|^p \, d\tau \\
   &\lesssim \mathbb{E} \int_0^t |\varepsilon^m |^p \,d\tau + \delta_m,
\end{align*}
wherein the mean value inequality to $x \mapsto \psi(x)x$ is used to show that
\[
| \psi(U^m)U^m - \psi(U)U | \le C(\psi) | U^m - U| = C|\varepsilon^m|,
\]
and 
\[ \delta_m = \mathbb{E} \int_0^t |\varphi_m(\eta^m) - 1|^p |\psi(U^m)|^p \, d\tau \rightarrow 0, \]
as $m\rightarrow\infty$, by the dominated convergence theorem.  Similarly, for the last term of the sum in (\ref{e:zs}),
\begin{align*}
\mathbb{E} \bigg| \int_0^t  &\left( \varphi_m^2(g^n)\psi^2(U^m) U^m -  \psi^2(U) U \right) 
     \operatorname{ad}^2_{X_i} \,d\tau \bigg|^p \\
   &\lesssim \mathbb{E} \int_0^t \left| \psi^2(U^m) U^m - \psi^2(U)U  \right|^p d\tau 
      + \mathbb{E} \int_0^t |\varphi_m^2(\eta^m) - 1|^p |\psi(U^m)|^{2p} \, d\tau \\
   &\lesssim \mathbb{E} \int_0^t |\varepsilon^m|^p \,d\tau + \delta_m,
\end{align*}
where the mean value inequality has now been applied to the function $x\mapsto\psi^2(x)x$, and 
\[ \delta_m = \mathbb{E} \int_0^t |\varphi_m^2(\eta^m) - 1|^p |\psi(U^m)|^{2p} \, d\tau \rightarrow 0, \]
as $m\rightarrow\infty$.  For the second term,
\begin{align*}
\mathbb{E} \bigg| \int_0^t  \left\langle \tilde{\nabla}\varphi_m(\eta^m), {X_i}\right\rangle \psi(U^m) U^m &\operatorname{ad}_{X_i}\,d\tau \bigg|^p \\
   &\lesssim \mathbb{E} \int_0^t \left| \left\langle \tilde{\nabla}\varphi_m(\eta^m), {X_i}\right\rangle \psi(U^m) U^m \right|^p \,d\tau \\
   &= \mathbb{E} \int_0^t \left| \left\langle \tilde{\nabla}\varphi_m(\eta^m), {X_i}\right\rangle \psi(U^m) (\varepsilon^m + U) \right|^p 
     \,d\tau \\
   &\lesssim \mathbb{E} \int_0^t |\varepsilon^m|^p \,d\tau + \delta_m,
\end{align*}
where
\[ \delta_m =  \mathbb{E} \int_0^t \left| \left\langle \tilde{\nabla}\varphi_m(\eta^m), {X_i}\right\rangle \psi(U^m) U \right|^p 
     \,d\tau \rightarrow 0, \]
as $m\rightarrow\infty$, since $\lim_{m\rightarrow\infty}|\tilde{\nabla}\varphi_m|=0$.
Finally, for the the third term, note first that
\begin{align*}
\varphi_m(g^n) \langle \psi&'(U^m) , U^m\operatorname{ad}_{X_i} \rangle U^m  \\
   &= \varphi_m(\eta^m) \left\langle \psi'(U^m) , U^m \operatorname{ad}_{X_i} \right\rangle (\varepsilon^m + U) \\
   &= \varphi_m(\eta^m) \left\langle \psi'(U^m) , U^m \operatorname{ad}_{X_i} \right\rangle\varepsilon^m 
       + \varphi_m(\eta^m) \left\langle \psi'(U^m) , (\varepsilon^m+U) \operatorname{ad}_{X_i} \right\rangle U \\
   &= \varphi_m(\eta^m) \left\langle \psi'(U^m) , U^m \operatorname{ad}_{X_i} \right\rangle\varepsilon^m 
       + \varphi_m(\eta^m) \left\langle \psi'(U^m) , \varepsilon^m \operatorname{ad}_{X_i} \right\rangle U \\
   & \qquad  + \varphi_m(\eta^m) \left\langle \psi'(U^m) , U \operatorname{ad}_{X_i} \right\rangle U.
\end{align*}
Thus,
\begin{align*}
\mathbb{E} \bigg| \int_0^t (\varphi_m&(\eta^m)\psi(U^m) \left\langle \psi'(U^m) , U^m\operatorname{ad}_{X_i} \right\rangle U^m 
     - \psi(U)\left\langle \psi'(U), U\operatorname{ad}_{X_i} \right\rangle U )\operatorname{ad}_{X_i}  \,d\tau \bigg|^p \\
   &\lesssim \mathbb{E} \int_0^t \left| \varphi_m(\eta^m) \left\langle \psi'(U^m) , U^m\operatorname{ad}_{X_i} \right\rangle\varepsilon^m 
      + \varphi_m(\eta^m) \left\langle \psi'(U^m) , \varepsilon^m \operatorname{ad}_{X_i} \right\rangle U \right.\\
   & \qquad \left. + \varphi_m(\eta^m) \left\langle \psi'(U^m) , U \operatorname{ad}_{X_i} \right\rangle U 
      - \left\langle \psi'(U), U\operatorname{ad}_{X_i} \right\rangle U \right|^p \,d\tau \\
   &\lesssim \mathbb{E} \int_0^t |\varepsilon^m|^p \,d\tau + \delta_m,
\end{align*}
where
\[ 
\delta_m = \mathbb{E} \int_0^t |\varphi_m(\eta^m) \left\langle \psi'(U^m) , U\operatorname{ad}_{X_i} \right\rangle - \left\langle \psi'(U), U\operatorname{ad}_{X_i} \right\rangle|^p |U|^p \,d\tau \rightarrow 0,
\]
as $m\rightarrow\infty$, since $\varphi_m(\eta^m)\psi'(U^m)\rightarrow\psi'(U)$ boundedly.  These bounds then imply that
\[ \mathbb{E} \sup_{\tau\le t}|\varepsilon_\tau^m|^p \le C\int_0^t |\varepsilon^m|^p\,d\tau + \delta_m, \]
for all $t\in[0,1]$.  Thus, by Gronwall's inequality, 
\[
\mathbb{E} \sup_{\tau\le1} |U^m_\tau - U_\tau|^p = \mathbb{E} \sup_{\tau\le1} |\varepsilon_\tau^m|^p 
   \le \delta_m e^{C} \rightarrow 0,
\]
as $m\rightarrow\infty$.
\end{proof}

\begin{prop}
\label{p:d2.5}
Let $\theta_t^m:\mathscr{W}\rightarrow \mathfrak{g}$ be as in Equation (\ref{e:d2.z}) with $u$ replaced by $u_m$; that is,
\begin{equation}
\label{e:thetan}
\theta_t^m 
    = \int_0^t \left( \left\langle \hat{\nabla}u_m(\eta^m) , \theta\right\rangle \operatorname{Ad}_{\eta^m} \circ d\vec{b} + u_m(\eta^m)\operatorname{Ad}_{\eta^m}
    d\vec{h}\right).
\end{equation}
Then
\[
\lim_{m\rightarrow\infty} \mathbb{E} \sup_{\tau\le1} | \theta^m_\tau - \theta_\tau|^p = 0,
\]
for all $p\in(1,\infty)$, where $\theta_t:\mathscr{W}\rightarrow\mathfrak{g}$ is the solution to 
\begin{equation}
\label{e:theta}
\theta_t = \int_0^t \left( \left\langle\hat{\nabla}v(\eta) , \theta \right\rangle \operatorname{Ad}_\eta  \circ \,d\vec{b} + v(\eta)\operatorname{Ad}_\eta \,d\vec{h} \right)
\end{equation}
\end{prop}

\begin{proof}
Let $U_t^m = \operatorname{Ad}_{\eta^m_t}$ and $U_t = \operatorname{Ad}_{\eta_t}$.  Rewrite Equation (\ref{e:thetan}) in It\^o form as
\begin{multline*}
d\theta^m = \left\langle\hat{\nabla} u_m(\eta^m) , \theta^m\right\rangle U^m d\vec{b} + u_m(\eta^m) U^m d\vec{h} \\
   + \frac{1}{2}\sum_{i=1}^k \left\langle \tilde{\nabla}\hat{\nabla}u_m\left(  \eta^m\right)  ,X_i\otimes
    \theta^m \right\rangle U^m X_i \,dt. 
\end{multline*}
Note that, formally, $\theta$ is the solution to Equation (\ref{e:d2.z}) with $u$ replaced by $v$ (although $v$ does not have compact support), and Equation (\ref{e:theta}) in It\^o form is
\begin{align*} 
d\theta &= \left\langle\hat{\nabla}v(\eta) , \theta \right\rangle \operatorname{Ad}_\eta  \,d\vec{b} + v(\eta)\operatorname{Ad}_\eta \,d\vec{h} 
      + \frac{1}{2}\sum_{i=1}^k \left\langle\tilde{\nabla}\hat{\nabla}v(\eta) , X_i\otimes\theta\right\rangle \operatorname{Ad}_\eta X_i \,dt, \\
   &= \left\langle\hat{\nabla}v(\eta) , \theta \right\rangle U  \,d\vec{b} + v(\eta)U \,d\vec{h} 
      + \frac{1}{2}\sum_{i=1}^k \left\langle\tilde{\nabla}\hat{\nabla}v(\eta) , X_i\otimes\theta\right\rangle U X_i \,dt.
\end{align*}

Let $\varepsilon^m := \theta^m - \theta$ (so that $\theta^m=\varepsilon^m+\theta$).  Then
\begin{multline*}
d\varepsilon^m = \left[ \left\langle\hat{\nabla} u_m(\eta^m),\theta^m\right\rangle U^m - \left\langle\hat{\nabla}v(\eta),\theta\right\rangle U\right] d\vec{b} 
      + \left[ u_m(\eta^m)U^m - v(\eta)U \right] d\vec{h} \\
   + \frac{1}{2}\sum_{i=1}^k \left[ \left\langle\tilde{\nabla}\hat{\nabla} u_m(\eta^m) , X_i\otimes\theta^m\right\rangle U^m 
   - \left\langle\tilde{\nabla}\hat{\nabla}v(\eta) , X_i\otimes\theta\right\rangle U\right] X_i \,dt.
%   &= \left[ \left\langle\hat{\nabla} u_m(g^n) , \theta + \varepsilon^m\right\rangle U^m-\left\langle\hat{\nabla}v(g) ,\theta\right\rangle U\right]  db
%    +\left[  u_m\left(  g^n\right)  U^m-v(g)  U\right]  dh \\
%  &\qquad + \frac{1}{2}\sum_{i=1}^k \left[ \left\langle\tilde{\nabla}\hat{\nabla} u_m(g^n) ,  %X_i\otimes(\theta+\varepsilon^m) \right\rangle U^m 
%    - \left\langle\tilde{\nabla}\hat{\nabla}v(g) , X_i\otimes\theta\right\rangle U\right] X_i \,dt.
\end{multline*}
Considering the first term of this expression, 
\begin{align*}
\left\langle\hat{\nabla} u_m(\eta^m),\theta^m\right\rangle& U^m -
    \left\langle\hat{\nabla}v(\eta),\theta\right\rangle U \\
    &=  \big\langle\hat{\nabla}u_m(\eta^m) , \theta + \varepsilon^m\big\rangle U^m 
      - \left\langle\hat{\nabla}v(\eta)  ,\theta\right\rangle U \\
    &= \left\langle\hat{\nabla}u_m(\eta^m) , \varepsilon^m\right\rangle U^m + \left\langle \hat{\nabla}u_m(\eta^m) , \theta\right\rangle U^m 
      - \left\langle\hat{\nabla}v(\eta) , \theta\right\rangle U.
%   &=   \left\langle\hat{\nabla}u_m\left(  \eta^m\right)  ,\varepsilon^m\right\rangle U^m\\ 
%   & \qquad  + \left\langle v(\eta^m)  \hat{\nabla}\varphi_m(\eta^m) , \theta\right\rangle U^m 
%     + \left\langle\varphi_m(\eta^m)\hat{\nabla}v(\eta^m)  ,\theta\right\rangle U^m-\left\langle\hat{\nabla}v(\eta) , \theta\right\rangle U.
\end{align*}
Using again that $U$ and $U^m$ remain in a fixed compact subset of $\mathrm{End}(\mathfrak{g})$, Proposition \ref{p:d2.55}, and the fact that $\hat{\nabla}u_m(\eta^m) \rightarrow \hat{\nabla}v(\eta)$ boundedly, 
\[
\mathbb{E} \left| \int_0^t\left[  \left\langle\hat{\nabla}u_m(\eta^m) , \theta^m \right\rangle U^m 
   - \left\langle\hat{\nabla}v(\eta) , \theta\right\rangle U \right] \, d\vec{b} \right|^p \lesssim \mathbb{E} \int_0^t |\varepsilon^m|^p \,d\tau + 
   \delta_m,
\]
where 
\[ \delta_m = \mathbb{E}\int_0^t \left|  \left\langle \hat{\nabla}u_m(\eta^m) , \theta\right\rangle U^m - \left\langle\hat{\nabla}v(\eta) , \theta\right\rangle U
   \right| d\tau \rightarrow 0, \]
as $m\rightarrow\infty$.  The second term converges to $0$ since
\[
\mathbb{E} \left| \int_0^t \left[ u_m(\eta^m) U^m - v(\eta)U \right] d\vec{h} \right|^p
  = \mathbb{E} \left| \int_0^t\left[  \varphi_m(\eta^m) v(\eta^m) U^m - v(\eta)U \right] d\vec{h} \right|^p \rightarrow 0,
\]
as $m\rightarrow\infty$, by the dominated convergence theorem.  For the third term, note that
\begin{multline*}
\left| \left[  \left\langle\tilde{\nabla}\hat{\nabla}u_m(\eta^m)  ,X_i\otimes (\theta+\varepsilon^m) \right\rangle U^m 
     - \left\langle\tilde{\nabla}\hat{\nabla}v(\eta)  ,X_i\otimes\theta\right\rangle U\right]  X_i\right| \\
   \lesssim |\varepsilon^m| + \left| \left[ \left\langle\tilde{\nabla}\hat{\nabla}u_m(\eta^m) , X_i\otimes\theta \right\rangle U^m 
     - \left\langle\tilde{\nabla}\hat{\nabla}v(\eta) , X_i\otimes\theta\right\rangle U\right]  X_i\right|,
\end{multline*}
and so
\begin{multline*}
\mathbb{E} \left| \int_0^t\left[  \left\langle\tilde{\nabla}\hat{\nabla} u_m(\eta^m) , X_i\otimes \theta^m
     \right\rangle U^m-\left\langle\tilde{\nabla}\hat{\nabla}v(\eta),X_i\otimes\theta\right\rangle U \right]  X_i \,d\tau\right|^p \\
   \lesssim \mathbb{E} \int_0^t |\varepsilon^m|^p \,d\tau + \delta_m,
\end{multline*}
where
\[
\delta_m = \mathbb{E} \left| \int_0^t\left[  \left\langle\tilde{\nabla}\hat{\nabla} u_m(\eta^m) , X_i \otimes\theta\right\rangle U^m
   - \left\langle\tilde{\nabla}\hat{\nabla}v(\eta) , X_i\otimes\theta\right\rangle U \right]  X_i \,d\tau\right|^p \rightarrow 0,
\]
as $m\rightarrow\infty$, by  the dominated convergence theorem, since $\tilde{\nabla}\hat{\nabla} u_m(\eta^m) \rightarrow \tilde{\nabla}\hat{\nabla}v(\eta)$ boundedly.

Putting these bounds together then shows
\[
\mathbb{E}\sup_{\tau\le t} |\varepsilon_\tau^m|^p \le C\mathbb{E} \int_0^t |\varepsilon^m|^p \,d\tau + \delta_m,
\]
for all $t\in[0,1]$, and again applying Gronwall's inequality gives
\[ \mathbb{E} \sup_{\tau\le1} |\theta_\tau^m - \theta_\tau |^p = \mathbb{E} \sup_{\tau\le1} |\varepsilon_\tau^m|^p 
   \le \delta_me^{C} \rightarrow 0, \]
as $m\rightarrow\infty$, finishes the proof.
\end{proof}

\begin{prop}
\label{p:Dh.eta}
Let $\eta^m$ be the solution to Equation (\ref{e:d2.7}) and $\theta^m$ be the solution to Equation (\ref{e:thetan}).  Then, for any $h\in \mathscr{H}$, $f\in C_c^\infty(G)$, and $t\in[0,1]$,
\[  \partial_h f(\eta_t^m) = \left\langle \hat{\nabla}f(\eta_t^m), \theta_t^m \right\rangle. \]
Furthermore, for $\eta$ the solution to (\ref{e:d2.8}) and $\theta$ the solution to (\ref{e:theta}),
\[ \lim_{m\rightarrow\infty} \mathbb{E} \left|\partial_h f(\eta^m_t) -  \left\langle \hat{\nabla}f(\eta_t), \theta_t \right\rangle\right|^p=0, \]
for all  $p\in(1,\infty)$.
\end{prop}

\begin{proof}
The first claim follows immediately from Proposition \ref{p:Dh}.  Now note that
\begin{align*}
\big|\partial_h f(\eta^m) - &\left\langle \hat{\nabla}f(\eta), \theta \right\rangle\big| 
   = \left|\left\langle \hat{\nabla}f(\eta^m), \theta^m \right\rangle - \left\langle \hat{\nabla}f(\eta), \theta \right\rangle\right| \\
   &\le \left|\left\langle \hat{\nabla}f(\eta^m), \theta^m \right\rangle - \left\langle \hat{\nabla}f(\eta^m), \theta \right\rangle\right| +
      \left|\left\langle \hat{\nabla}f(\eta^m), \theta \right\rangle - \left\langle \hat{\nabla}f(\eta), \theta \right\rangle\right| \\
   &\le  |\hat{\nabla}f|| \theta^m - \theta | + |\hat{\nabla}f(\eta^m)- \hat{\nabla}f(\eta)||\theta|.
\end{align*}
Thus,
\begin{multline*} 
\lim_{m\rightarrow\infty} \mathbb{E} \left| \partial_h f(\eta^m) - \left\langle \hat{\nabla}f(\eta), \theta \right\rangle \right|^p \\
   \le \lim_{m\rightarrow\infty} \mathbb{E}\left[|\hat{\nabla}f|| \theta^m - \theta | + |\hat{\nabla}f(\eta^m) -
      \hat{\nabla}f(\eta)||\theta|\right]^p = 0, 
\end{multline*}
by Propositions \ref{p:eta} and \ref{p:d2.5} and the dominated convergence theorem.
\end{proof}

\begin{cor}
\label{c.d2.6}
For any $h\in \mathscr{H}$, $f\in C_c^\infty(G)$, and $t\in[0,1]$, $f(\eta_t)\in \mathrm{Dom}(\partial_h)$ and 
\[ \partial_h f(\eta_t) = \left\langle \hat{\nabla}f(\eta_t), \theta_t \right\rangle\in L^{\infty-}(\mu). \]
\end{cor}

This corollary follows from $\partial_h$ being a closed operator (taking $\partial_h=\overline{\partial}_h$).
%Propositions \ref{p:eta} and \ref{p:Dh.eta} together imply that $f(\eta_t)\in\mathrm{Dom}(\partial_h)$, for any %$f\in C_c^\infty(G)$, and 
%\[ \partial_h f(\eta_t) = \left\langle \hat{\nabla}f(\eta_t), \theta \right\rangle \in L^p(\mu), \]
%for all $p\in[1,\infty)$.  
%By Theorem \ref{t:diff}, this implies that $f(\eta_t)\in\mathcal{D}^{1,p}$, for all $p\in(1,\infty)$.
Now removing the cutoff functions completes the proof of the primary result of this paper.  

\iffalse
\begin{thm}
\label{t:xismooth}
Let $\xi:[0,1]\times \mathscr{W}\rightarrow G$ denote the solution to Equation (\ref{e:SDE}), and let $\Theta:\mathscr{W}\rightarrow\mathfrak{g}$ be the solution to
\begin{equation}
\label{e:Theta}
\Theta_t := \int_0^t \operatorname{Ad}_{\xi} d\vec{h} = \int_0^t \operatorname{Ad}_{\xi_\tau} X_i\dot{h}_\tau^i \,d\tau.
\end{equation}
Then, for any $f\in C_c^\infty(G)$ and $t\in[0,1]$, $f(\xi_t)\in\mathcal{D}^{1,\infty}$, and for any $h\in \mathscr{H}$,
\begin{equation}
\label{e:xi'}
\partial_h [f(\xi_t)] = \int_0^t \left\langle \hat{\nabla}f(\xi_t), \operatorname{Ad}_{\xi_\tau} X_i \right\rangle \dot{h}^i_\tau \,d\tau.
\end{equation}
Thus,
\[ (D[f(\xi_t)])^i = \int_0^{\cdot\wedge t} \left\langle \hat{\nabla}f(\xi_t), \operatorname{Ad}_{\xi_\tau} X_i \right\rangle\,d\tau, \]
componentwise in $\mathscr{H}$.
\end{thm}
\fi

\hspace*{0.01in}
{\bf Proof of Theorem \ref{t:xi.smooth}.} 
Let $\psi\in C_c^\infty(\mathrm{End}(\mathfrak{g}),[0,1])$ be as in Notation \ref{n:phim}, and define $v_n(g) := \psi(n^{-1}\operatorname{Ad}_g)$.  Let $\xi^n:[0,1]\times \mathscr{W}\rightarrow G$ be the solution to the Stratonovich equation
\[
d\xi^n=v_n(\xi^n)\xi^n \circ d\vec{b},  \text{ with } \xi^n_0 = e.
\]
By Lemma \ref{l.d2.1}, $\xi^n$ exists for all time $t\in[0,1]$.  Noting that $v_n\rightarrow1$ and $\tilde{\nabla}v_n\rightarrow0$ boundedly as $n\rightarrow\infty$, an argument identical to that in Proposition \ref{p:eta} shows that, for all $f\in C_c^\infty(G)$,
\begin{equation}
\label{e:xint}
\lim_{n\rightarrow\infty} \mathbb{E}\sup_{\tau\le1} |f(\xi_\tau^n)-f(\xi_\tau)|^p = 0,
\end{equation}
for all $p\in(1,\infty)$.  Note that this convergence implies that $f(\xi_t)\in L^{\infty-}(\mu)$.

Now let $\Theta:\mathscr{W}\rightarrow\mathfrak{g}$ denote the solution to
\begin{equation}
\label{e:Theta}
\Theta_t := \int_0^t \operatorname{Ad}_{\xi} d\vec{h} = \int_0^t \operatorname{Ad}_{\xi_\tau} X_i\dot{h}_\tau^i \,d\tau,
\end{equation}
and let $\Theta_t^n:\mathscr{W}\rightarrow\mathfrak{g}$ denote the solution to the It\^o equation
\begin{multline}
\label{e:Thetan}
d\Theta^n  =\left\langle\hat{\nabla}v_n\left(  \xi^n\right)  ,\Theta^n\right\rangle \operatorname{Ad}_{\xi^n} d\vec{b} + v_n(\xi^n) \operatorname{Ad}_{\xi^n} \,d\vec{h} \\
     +\frac{1}{2} \sum_{i=1}^k \left\langle \tilde{\nabla}\hat{\nabla}v_n(\xi^n) , X_i\otimes\Theta^n \right\rangle \operatorname{Ad}_{\xi^n} X_i \,dt,
\end{multline}
with $\Theta_0^n=0$.  Now it must be shown that
\begin{equation}
\label{e:n}
\lim_{n\rightarrow\infty}\mathbb{E} \sup_{\tau\le1} |\Theta^n_\tau - \Theta_\tau|^p = 0,
\end{equation}
for all $p\in(1,\infty)$.  So let $W^n=\operatorname{Ad}_{\xi^n}$ and $W=\operatorname{Ad}_\xi$.  Then $W^n,W:\mathscr{W}\rightarrow\mathrm{End}(\mathfrak{g})$ satisfy the equations
\begin{equation}
\label{e:Wn}
dW^n = v_n(\xi^n) W^n\circ\operatorname{ad}_{d\vec{b}}, \text{ with } W_{0}^n=I,
\end{equation}
and
\begin{equation}
\label{e:W}
dW = W \circ\operatorname{ad}_{d\vec{b}}, \text{ with } W_{0}=I.
\end{equation}
Rewriting Equation (\ref{e:Thetan}), $\Theta^n$ solves
\begin{multline*}
d\Theta^n  =\left\langle\hat{\nabla}v_n\left(  \xi^n\right)  ,\Theta^n\right\rangle W^n d\vec{b} + v_n(\xi^n) W^n \,d\vec{h} \\
     +\frac{1}{2} \sum_{i=1}^k \left\langle \tilde{\nabla}\hat{\nabla}v_n(\xi^n) , X_i\otimes\Theta^n \right\rangle W^n X_i \,dt,
\end{multline*}
and Equation (\ref{e:Theta}) implies that $d\Theta = W \,d\vec{h}$.
Thus, for $\varepsilon^n := \Theta^n-\Theta$, 
\begin{multline}
\label{e:ggg}
d\varepsilon^n = \left\langle\hat{\nabla}v_n\left(  \xi^n\right)  ,\Theta^n\right\rangle W^n d\vec{b} + (v_n(\xi^n) W^n - W) d\vec{h} \\
     + \frac{1}{2} \sum_{i=1}^k  \left\langle \tilde{\nabla}\hat{\nabla}v_n(\xi^n) , X_i\otimes\Theta^n \right\rangle W^n X_i \,dt.
\end{multline}
Since $W^n$ and $W$ solve Equations (\ref{e:Wn}) and (\ref{e:W})
which have smooth, bounded coefficients, Theorem V-10.1 of Ikeda and Watanabe \cite{IkedaWatanabe89} implies that $W$ and $W^n$ are in $\mathcal{D}^\infty(\mathrm{End}(\mathfrak{g}))$.  Standard matrix group results imply that
\begin{equation*}
\lim_{n\rightarrow\infty} \mathbb{E} \sup_{\tau\le1} |W^n_\tau - W_\tau|^p 
\qquad \text{and}\qquad
\lim_{n\rightarrow\infty} \mathbb{E} \sup_{\tau\le1} |\partial_h W^n_\tau - \partial_h W_\tau|^p,
\end{equation*}
for all $p\in(1,\infty)$; see for example Propositions 6 and 8 in \cite{app}.  Furthermore,
\[
\partial_h W^n=\partial_h \operatorname{Ad}_{\xi^n} = \operatorname{ad}_{\Theta^n}\operatorname{Ad}_{\xi^n} = \operatorname{ad}_{\Theta^n} W^n.
\]
This then implies that
\begin{align*}
\left\langle\hat{\nabla}v_n(\xi^n),\Theta^n \right\rangle &= \frac{d}{dt}\bigg|_0 \psi(n^{-1}\operatorname{Ad}_{e^{t\Theta^n}\xi^n})  \\
   &= \frac{1}{n} \left\langle \psi'(n^{-1}\operatorname{Ad}_{\xi^n}) , \operatorname{ad}_{\Theta^n}\operatorname{Ad}_{\xi^n}\right\rangle  
   =\frac{1}{n} \left\langle \psi'(n^{-1}W^n) , \partial_h W^n\right\rangle  .
\end{align*}
Thus, for the first term of Equation (\ref{e:ggg}),
\begin{align*}
\mathbb{E} \bigg| \int_0^t \big\langle\hat{\nabla} v_n(\xi^n) &, \Theta^n\big\rangle W^n \,d\vec{b}\bigg|^p 
   \lesssim \mathbb{E} \int_0^t \frac{1}{n} \left| \left\langle \psi'(n^{-1}W^n) , \partial_h W^n \right\rangle  W^n\right|^p \,d\tau \\
   &\lesssim \mathbb{E} \int_0^t \left( \frac{1}{n} \left| \psi'(n^{-1}W^n)  \right| \left| W^n \right|\right)^p |\partial_h W^n|^p \,d\tau \\
   &\lesssim \mathbb{E} \int_0^t \left( \frac{1}{n} \left| \psi'(n^{-1}W^n)  \right| \left| W^n \right| \right)^p |\partial_h W|^p \,d\tau \\
   & \qquad + \mathbb{E} \int_0^t \left( \frac{1}{n} \left| \psi'(n^{-1}W^n)  \right| \left| W^n \right| \right)^p |\partial_h W^n - \partial_h W|^p \,d\tau
   \rightarrow 0,
\end{align*}
as $n\rightarrow\infty$, by the dominated convergence theorem.  Similarly,
\begin{align*}
\left\langle\tilde{\nabla}\hat{\nabla} v_n(\xi^n)  , X_i  \otimes \Theta^n\right\rangle 
   &= \frac{1}{n^2}\left\langle \psi''(n^{-1}W^n) , W^n\operatorname{ad}_{X_i }\otimes \operatorname{ad}_{\Theta^n}\operatorname{Ad}_{\xi^n} \right\rangle  \\
%   &= \frac{1}{n^2}\left\langle \psi''(n^{-1}W^n) , W^n\operatorname{ad}_{A}\otimes \operatorname{ad}_{\Theta^n}W^n\right\rangle  \\
   &= \frac{1}{n^2}\left\langle \psi''(n^{-1}W^n) , W^n\operatorname{ad}_{X_i }\otimes \partial_h W^n\right\rangle ,
\end{align*}
so that 
\begin{align*}
\mathbb{E} \bigg| \int_0^t  
      \big\langle\tilde{\nabla}\hat{\nabla} v_n&(\xi^n),  X_i \otimes \Theta^n \big\rangle W^nX_i  \,d\tau\bigg|^p \\
   &= \mathbb{E} \bigg| \int_0^t  \frac{1}{n^2}\left\langle \psi''(n^{-1}W^n) , W^n \operatorname{ad}_{X_i}\otimes 
      \partial_h W^n \right\rangle W^n X_i  \bigg|^p \,d\tau \\
   &\lesssim \mathbb{E} \int_0^t \left( \frac{1}{n^2} \left|\psi''(n^{-1}W^n)\right| |W^n|^2 \right)^p |\partial_h W^n|^p \,d\tau \\
   &\lesssim \mathbb{E} \int_0^t \left( \frac{1}{n^2} \left| \psi''(n^{-1}W^n)\right| |W^n|^2\right)^p |\partial_h W|^p \,d\tau \\
   &\qquad + \mathbb{E} \int_0^t \left( \frac{1}{n^2} \left|\psi''(n^{-1}W^n)\right| |W^n|^2 \right)^p |\partial_h W^n - \partial_h W|^p \,d\tau
   \rightarrow 0,
\end{align*}
as $n\rightarrow\infty$, again by dominated convergence.  Finally,
\begin{multline*}
\mathbb{E} \bigg| \int_0^t \left(v_n(\xi^n)W^n- W \right) \,d\vec{h} \bigg|^p \\
   \lesssim \mathbb{E} \int_0^t |v_n(\xi^n)-1|^p|W|^p \,d\vec{h} + \mathbb{E} \int_0^t \left|v_n(\xi^n)\right|^p\left|W^n- W\right|^p 
   \,d\vec{h}
   \rightarrow 0
\end{multline*}
as $n\rightarrow\infty$, by the dominated convergence theorem.  Thus, (\ref{e:n}) is verified.

Now, by Corollary \ref{c.d2.6}, for any $f\in C_c^\infty(G)$ and $t\in[0,1]$, $\partial_h[ f(\xi_t^n)]\in L^{\infty-}(\mu)$, and
\[ \partial_h [f(\xi^n_t)] = \left\langle \hat{\nabla}f(\xi^n_t),\Theta^n_t \right\rangle. \]
Thus, by the same argument as in Proposition \ref{p:Dh.eta}, Equations (\ref{e:xint}) and (\ref{e:n}) imply that
\[ \lim_{n\rightarrow\infty} \mathbb{E}\left|\partial_h f(\xi^n_t) -  \left\langle \hat{\nabla}f(\xi_t),\Theta_t \right\rangle \right|^p = 0, \]
for all $p\in(1,\infty)$.  Since $\partial_h$ is a closed operator, this and $f(\xi_t)\in L^{\infty-}(\mu)$ together imply that $f(\xi_t)\in\mathrm{Dom}(\partial_h)$ and 
\[ \partial_h [f(\xi_t)] = \left\langle \hat{\nabla}f(\xi_t), \Theta_t \right\rangle =  \left\langle \hat{\nabla}f(\xi_t), 
   \int_0^t \operatorname{Ad}_{\xi_\tau} X_i\dot{h}_\tau^i \,d\tau \right\rangle \in L^{\infty-}(\mu). \]
In particular, for any $h\in \mathscr{H}$ such that $\|h\|_\mathscr{H}=1$,
\[ \mathbb{E}|\partial_h [f(\xi_t)]|^p \le \|\hat{\nabla}f \|_\infty \mathbb{E}\left| \int_0^t \operatorname{Ad}_{\xi_\tau} X_i\dot{h}_\tau^i \,d\tau \right|^p
   \le \|\hat{\nabla}f \|_\infty \sum_{i=1}^k \mathbb{E}\int_0^t \left|\operatorname{Ad}_{\xi_\tau} X_i \right|^p \,d\tau. \] 
Since
\[
\mathbb{E} \sup_{\tau\le1} |\operatorname{Ad}_{\xi_\tau}|^p \le C, 
\]
for some finite constant $C$ (see, for example, the proof of Proposition 3 in \cite{app}),  $\|f(\xi_t)\|_{G^{1,p}}<\infty$ for all $p\in(1,\infty)$, where 
\[ \|F\|_{G^{1,p}} = \|F\|_{L^p(\mu)} + \left(\mathbb{E} \sup_{\|h\|_\mathscr{H}=1} |\partial_h F|^p \right)^{1/p} ; \] 
see Section \ref{s:Mall.def}.  Then by Theorem \ref{t:diff}, $f(\xi_t)\in\mathcal{D}^{1,\infty}$, and
\[
\mathbb{E}[f(\xi_t)D^*h] %= \mathbb{E}[(D[f(\xi_t)],h)_H]
   = \mathbb{E}[\partial_h [f(\xi_t)]] 
%   = \mathbb{E}\left[\left\langle \hat{\nabla}f(\xi_t), \int_0^t \operatorname{Ad}_{\xi_\tau} \,d\vec{h}_\tau \right\rangle\right]
   = \mathbb{E} \left[\int_0^t \left\langle \hat{\nabla}f(\xi_t), \operatorname{Ad}_{\xi_\tau} X_i \right\rangle \dot{h}_\tau^i \,d\tau\right],
\]
%and so $f(\xi_t)\in \mathcal{D}(D^{**})= \mathcal{D}(D)$, since $D$ closed implies that $D^{**}=D$ (where we are %thinking $D=\bar{D}$), 
implies that
\begin{equation}
\label{e:D} 
(D_s[f(\xi_t)])^i =\left\langle \hat{\nabla}f(\xi_t),  \int_0^{s\wedge t} \operatorname{Ad}_{\xi_\tau} X_i \,d\tau \right\rangle, 
\end{equation}
componentwise in $\mathscr{H}$. 

Finally, $\overline{W} = \int_0^\cdot \operatorname{Ad}_{\xi_\tau} d\tau \in \mathcal{D}^\infty(\mathscr{H}(\mathrm{End}(\mathfrak{g})))$ (see, for example, Proposition 5 in \cite{app}).  Combined with Equation (\ref{e:D}), this shows that $f(\xi_t)\in\mathcal{D}^\infty$, for all $f\in C_c^\infty(G)$ and $t\in[0,1]$.
\hfill{\tiny $\blacksquare$}

\bibliographystyle{amsplain}
\bibliography{biblio}

\end{document}